\documentclass[11pt]{amsart}
\usepackage{amssymb,amsmath,amsthm,amsfonts,mathrsfs,mathtools}
\usepackage[hmargin=3cm,vmargin=3.5cm]{geometry}
\usepackage{xcolor}
\definecolor{darkblue}{rgb}{0,0,0.4} 
\usepackage[colorlinks=true, citecolor=darkblue, filecolor=darkblue, linkcolor=darkblue,urlcolor=darkblue]{hyperref}
\usepackage{graphicx,psfrag}
\usepackage{amscd}  
\usepackage[shellescape]{gmp}
\usepackage{stmaryrd}  
\usepackage[dvips]{epsfig}
\usepackage{tikz}
\usepackage[textsize=tiny]{todonotes}

\newtheorem{theorem}{Theorem}[section]

\newtheorem{remark}[theorem]{Remark}

\newtheorem{corollary}{Corollary}

\renewcommand{\th}{^{\text{th}}}

\let\oldtocsection=\tocsection
\let\oldtocsubsection=\tocsubsection
\renewcommand{\tocsection}[2]{\hspace{0em}\oldtocsection{#1}{#2}}
\renewcommand{\tocsubsection}[2]{\hspace{1em}\oldtocsubsection{#1}{#2}}

\usepackage{chngcntr}
\counterwithin{figure}{section}

\title{Categorical lifting of the Jones polynomial: a survey}

\author{Mikhail Khovanov}
 \address{Department of Mathematics, Columbia University, New York, NY 10027, USA}
 \email{\href{mailto:khovanov@math.columbia.edu}{khovanov@math.columbia.edu}}
 \author{Robert Lipshitz}
 \address{Department of Mathematics, University of Oregon,  Eugene, OR 97403}
 \email{\href{mailto:lipshitz@uoregon.edu}{lipshitz@uoregon.edu}}
 
\date{September 20, 2021}

\begin{document}

\def\R{\mathbb R}
\def\Q{\mathbb Q}
\def\Z{\mathbb Z}
\def\N{\mathbb N} 
\def\C{\mathbb C}
\def\S{\mathbb S}
\def\CP{\mathbb P}
\def\FF{\mathbb F}
\renewcommand\SS{\ensuremath{\mathbb{S}}}
\newcommand{\ophana}{\overline{\phantom{a}}}
\newcommand{\SU}{\mathit{SU}}

\def\mathcenter#1{%
  \vcenter{\hbox{$#1$}}%
}

\def\l{\lbrace}
\def\r{\rbrace}
\def\o{\otimes}
\def\lra{\longrightarrow}
\def\Ext{\mathrm{Ext}}
\def\mc{\mathcal}
\def\mf{\mathfrak} 
\def\mcC{\mathcal{C}}
\def\mcI{\mathcal{I}}
\def\uFr{\underline{\mathrm{Fr}}}

\def\bbb{\mathbb{B}}
\def\ovb{\overline{b}}
\def\tr{{\sf tr}} 
\def\det{{\sf det }} 

\def\lra{\longrightarrow}
\def\kk{\mathbf{k}}  
\def\gdim{\mathrm{gdim}}  
\def\rk{\mathrm{rk}}

\def\Cob{\mathrm{Cob}} 

\newcommand{\germ}{\mathfrak}
\newcommand{\Kh}{\mathrm{H}}
\newcommand{\rKh}{\widetilde{\Kh}}
\newcommand{\HF}{\mathit{HF}}
\newcommand{\HFa}{\widehat{\HF}}
\newcommand{\CF}{\mathit{CF}}
\newcommand{\CFa}{\widehat{\CF}}
\newcommand{\brak}[1]{\langle #1 \rangle}

\newcommand{\rank}{\mathrm{rank}}

\newcommand{\oplusop}[1]{{\mathop{\oplus}\limits_{#1}}}
\newcommand{\ang}[1]{\langle #1 \rangle } 

\newcommand{\mcA}{{\mathcal A}}
\newcommand{\mcR}{{\mathcal R}}

\newcommand{\undn}{\mathbf{n}}
\newcommand{\Hom}{\mbox{Hom}}
\newcommand{\id}{\mbox{id}}
\newcommand{\Id}{\mbox{Id}}
\newcommand{\End}{\mathrm{End}}
\newcommand{\Kaum}{\mathrm{Kau}} 
\newcommand{\co}{\nobreak\mskip2mu\mathpunct{}\nonscript
  \mkern-\thinmuskip{:}\penalty300\mskip6muplus1mu\relax}

\begin{abstract} This is a brief review of the categorification of the Jones polynomial and its significance and ramifications in geometry, algebra, and low-dimensional topology.
\end{abstract}

\dedicatory{This paper is dedicated to the memory of Vaughan Jones, whose insights have illuminated so many beautiful mathematical paths.}


\maketitle

%
%

\section{Constructions of the Jones polynomial}\label{sec:Jones}

The spectacular discovery by Vaughan Jones~\cite{Jones85:Bulletin,Jones87:Hecke} of  the Jones polynomial of links has led to many follow-up  developments in mathematics. In this note we will survey one of these developments, the discovery of a combinatorially-defined homology theory of links, functorial under link cobordisms in 4-space, and its connections to algebraic geometry, symplectic geometry, gauge theory, representation theory, and stable homotopy theory. 

The Jones polynomial $J(L)$ of an oriented link $L$ in $\R^3$ is determined uniquely by the skein relation 
%
\begin{equation}\label{eq_skein_1}
  q^{-2}J\Biggl(
  \mathcenter{
  \begin{tikzpicture}[scale=.8]
    \draw[dashed] (0,0) circle (1);
    \draw[->, thick] (-0.707,-.707) -- (.707,.707);
    \draw[thick] (0.707,-.707) -- (.1,-.1);
    \draw[->, thick] (-.1,.1) -- (-.707,.707);
  \end{tikzpicture}}
  \Biggr)
  -
  q^{2} J\Biggl(
  \mathcenter{
  \begin{tikzpicture}[scale=.8]
    \draw[dashed] (0,0) circle (1);
    \draw[->, thick] (0.707,-.707) -- (-.707,.707);
    \draw[thick] (-0.707,-.707) -- (-.1,-.1);
    \draw[->, thick] (.1,.1) -- (.707,.707);
  \end{tikzpicture}}
\Biggr)
=
(q^{-1}-q)J\Biggl(
  \mathcenter{
  \begin{tikzpicture}[scale=.8]
    \draw[dashed] (0,0) circle (1);
    \draw[->, thick] (-0.707,-.707) .. controls (-.1,0) .. (-.707,.707);
    \draw[->, thick] (0.707,-.707) .. controls (.1,0) .. (.707,.707);
  \end{tikzpicture}}
\Biggr)
\end{equation}
and the normalization that the polynomial of the unknot satisfies $J(U)=1$.  The multiplicativity property  $J(L\sqcup U)=(q+q^{-1})J(L)$, that is, that the disjoint union with the unknot scales  the  invariant by $q+q^{-1}$, suggests another natural normalization, $J(U)=q+q^{-1}$ and  $J(\emptyset)=1$, where $\emptyset$ is the empty link.

The polynomial $J(L)$ originally arose from Jones's work on $C^*$-algebras, where the braid relations and Temperley-Lieb relations appeared organically~\cite{Jones83:braid,Jones83:subfactors}. As we will see below, it also has connections to many other areas, from representation theory to gauge theory. Many of these connections first appeared or were foreshadowed in papers of Jones's, including the connections to quantum groups and statistical mechanics~\cite{Jones89:statistical}, Hecke algebras and traces~\cite{Jones83:braid,Jones87:Hecke}, and many other topics~\cite{Jones91:subfactors-book}. In addition to inspiring at least half a dozen different fields in mathematics, the Jones polynomial and its descendants have had remarkable applications to topology. Some we will touch on below; others, like its central role in resolving the famous Tait conjectures or its deep connections to hyperbolic geometry, we leave to other authors.

\begin{figure}
  \centering
  \[
    \left\langle
     \mathcenter{\begin{tikzpicture}[scale=.8]
        \draw[dashed] (0,0) circle (1);
        \draw[thick] (-.707,-.707) -- (.707,.707);
        \draw[thick] (.707,-.707) -- (.1,-.1);
        \draw[thick] (-.1,.1) -- (-.707,.707);
      \end{tikzpicture}}\right\rangle\quad=\quad
      q^{-1/2}\left\langle \mathcenter{\begin{tikzpicture}[scale=.8]
    \draw[dashed] (0,0) circle (1);
    \draw[-, thick] (-.707,-.707) .. controls (-.1,0) .. (-.707,.707);
    \draw[-, thick] (.707,-.707) .. controls (.1,0) .. (.707,.707);
  \end{tikzpicture}}\right\rangle\quad-\quad
  q^{1/2}\left\langle\mathcenter{\begin{tikzpicture}[scale=.8]
    \draw[dashed] (0,0) circle (1);
    \draw[-, thick] (-.707,-.707) .. controls (0,-.1) .. (.707,-.707);
    \draw[-, thick] (-.707,.707) .. controls (0,.1) .. (.707,.707);
  \end{tikzpicture}}\right\rangle
  \]
  \caption{\textbf{Kauffman bracket skein relation.} Given a diagram $D$ and a crossing in it as on the left, there are two ways $D_0$ and $D_1$ to resolve the crossing, as on the right. The Kauffman bracket of $D$, $D_0$, and $D_1$ are related as shown.}
  \label{fig:resolution}
\end{figure}
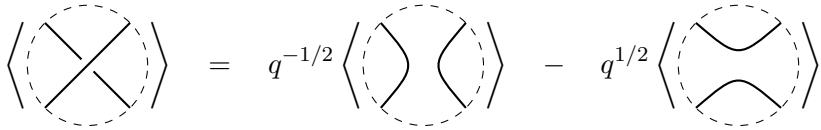

While it is fairly easy to see that at most one knot invariant satisfies Relation~(\ref{eq_skein_1}) and any given normalization for $J(U)$,  it is not immediately obvious that~(\ref{eq_skein_1}) is consistent. A simple way to see the existence of a  knot invariant satisfying~(\ref{eq_skein_1}) was discovered by L.~Kauffman~\cite{Kauffman87:state}.
Pick a planar diagram $D$ of $L$, forget about the orientation of $L$, and  
resolve each crossing of $D$ into a linear combination of two crossingless diagrams as shown in Figure~\ref{fig:resolution}. Any time a simple closed curve without crossings arises, remove it and scale the remaining term by $q+q^{-1}$.
The end result is a Laurent polynomial $\brak{D}\in q^{n/2}\Z[q,q^{-1}]$ (where $n$ is the number of crossings of $D$), the \emph{Kauffman bracket of $D$}. We can now bring back the orientation  of $L$ and scale $\brak{D}$ by a monomial in terms of the number $n_+$ of positive crossings and $n_-$ of negative crossings (the first and second pictures in Formula~\eqref{eq_skein_1}): 
\begin{equation}\label{eq:writhe-renorm}
    K(D)  \coloneqq (-1)^{n_-}q^{3(n_+-n_-)/2} \brak{D} \ \in \ \Z[q,q^{-1}]. 
\end{equation}
It is straightforward  to check that $K(D)$ is invariant under Reidemeister moves of oriented link diagrams, hence gives rise to a link invariant $K(L)$. Further, by applying the unoriented skein relation from Figure~\ref{fig:resolution} at the crossing of the two diagrams on the left of Relation~(\ref{eq_skein_1}), one sees that $K(L)$ satisfies Relation~(\ref{eq_skein_1}). So, we have:
\begin{theorem}
  (Kauffman \cite{Kauffman87:state}) For any oriented link $L$, $J(L)=K(L)$.
\end{theorem}

\section{Categorification of the Jones polynomial for links and tangles}\label{sec:Kh}
\subsection{Categorification for links}
E.~Witten showed~\cite{Witten89:Jones} at a physical level of rigor that the Chern-Simons path integral, with gauge group $\SU(2)$ and parameter $q$ a root of unity, gives rise to an invariant of 3-manifolds intricately  related to the  Jones polynomial. The case of gauge group $U(1)$ was considered earlier by A.~Schwarz, who showed that the path integral evaluates to the Reidemeister torsion~\cite{Schwarz77:CS}.
Shortly afterwards, N.~Reshetikhin and V.~Turaev~\cite{RT91:3mfld} gave a mathematically precise proof that suitable linear combinations of the Jones polynomial of cables of a framed link $L$, evaluated at $q$ an $N$-th root of unity, give invariants $\tau_N(M)$ of an oriented 3-manifold $M$ obtained by surgery on $L$; the resulting invariants are called \emph{Witten-Reshetikhin-Turaev invariants}.

Motivated by these developments and by constructions in geometric representation theory, notably by the work of G.~Lusztig~\cite{Lusztig90:bases} and A.~Beilinson, Lusztig, and R.~MacPherson~\cite{BLM90:qdeform}, L.~Crane and I.~B.~Frenkel conjectured~\cite{CF94:4d} that the  Witten-Reshetikhin-Turaev three-manifold invariant lifts to a 4-dimensional TQFT. They coined the  term  \emph{categorification} to  describe such a lifting of an $(n-1)$-dimensional TQFT to an $n$-dimensional TQFT. 

Despite many insights  into the possible structure of such a theory since then, its existence still  remains a conjecture. Nonetheless, the Crane-Frenkel  conjecture  motivated  the discovery of a categorification of the Jones polynomial by the first author~\cite{Khovanov00:Jones}.
In that  categorification,  the parameter $q$ becomes a grading shift of the \emph{quantum grading}, and the theory assigns to an oriented link $L\subset\R^3$ bigraded homology groups 
\begin{equation}
  \Kh(L) =\oplusop{i,j\in \Z} \Kh^{i,j}(L),
\end{equation}
functorial under smooth link cobordisms, and with the Jones polynomial as their Euler characteristic:
\begin{equation}
    J(L)=\sum_{i,j\in \Z}(-1)^iq^j \rank(\Kh^{i,j}(L)).
\end{equation}

A way to construct this theory can be guessed by lifting the Kauffman skein relation to a long exact sequence for homology. That is, up to appropriate grading shifts, there is an exact sequence 
\[
  \cdots\longrightarrow
  \Kh\left(\mathcenter{\begin{tikzpicture}[scale=.7]
        \draw[dashed] (0,0) circle (1);
        \draw[-, thick] (-.707,-.707) .. controls (0,-.1) .. (.707,-.707);
        \draw[-, thick] (-.707,.707) .. controls (0,.1) .. (.707,.707);
      \end{tikzpicture}}\right)\longrightarrow
    \Kh\left(
     \mathcenter{\begin{tikzpicture}[scale=.7]
        \draw[dashed] (0,0) circle (1);
        \draw[thick] (-.707,-.707) -- (.707,.707);
        \draw[thick] (.707,-.707) -- (.1,-.1);
        \draw[thick] (-.1,.1) -- (-.707,.707);
      \end{tikzpicture}}\right)\longrightarrow
      \Kh\left( \mathcenter{\begin{tikzpicture}[scale=.7]
    \draw[dashed] (0,0) circle (1);
    \draw[-, thick] (-.707,-.707) .. controls (-.1,0) .. (-.707,.707);
    \draw[-, thick] (.707,-.707) .. controls (.1,0) .. (.707,.707);
  \end{tikzpicture}}\right)\longrightarrow
\Kh\left(\mathcenter{\begin{tikzpicture}[scale=.7]
    \draw[dashed] (0,0) circle (1);
    \draw[-, thick] (-.707,-.707) .. controls (0,-.1) .. (.707,-.707);
    \draw[-, thick] (-.707,.707) .. controls (0,.1) .. (.707,.707);
  \end{tikzpicture}}\right)[1]\longrightarrow\cdots.
\]

Suppose further that, given a diagram $D$ for $L$, there is a chain complex $C(D)$ computing $\Kh(L)$ and the long exact sequence is induced by an isomorphism between the complex $C(D)$ and the cone of a map between $C(D_0)$ and $C(D_1)$ (where $D_0$ and $D_1$ are as in Figure~\ref{fig:resolution}). The Jones invariant of the unknot is $q+q^{-1}$, which is the graded rank of a free graded abelian group $A$ with generators in degrees $-1$ and $1$. The philosophy of topological quantum field theories then suggests to associate $A^{\otimes k}$ to a $k$-component unlink diagram. Natural maps $C(D_0)\to C(D_1)$ between these complexes for 
resolutions of $D$ can be obtained from a commutative Frobenius algebra structure on $A$: change of resolution is a cobordism, and Frobenius algebras correspond to $2$-dimensional topological quantum field theories, assigning maps to cobordisms between $1$-manifolds. It turns out that $A$ is unique up to obvious symmetries: with generators in (quantum) degrees $-1$ and $1$ denoted by $1$ and $X$, respectively, the multiplication $m$ and the trace $\epsilon$ on $A$ are given by 
\begin{align} 
    &   A = \Z  1 \oplus  \Z  X, \quad  1\cdot a = a\cdot 1= a \ (\forall a\in A), \\ 
    &   X\cdot X = 0, \qquad \epsilon(1)=0,\quad \epsilon(X)=1. 
\end{align}
Dualizing the  multiplication via $\epsilon$ leads to a comultiplication, with 
\begin{equation}
    \Delta(1)=1\otimes X + X\otimes 1, \qquad \Delta(X)=X\otimes X. 
\end{equation}
Explicitly, $m$ and $\Delta$ allow one to write down maps associated to all local topology changes between $2^n$ full resolutions of an $n$-crossing diagram $D$, giving a commutative $n$-dimensional cube with powers of $A$ at its vertices and maps $m$ and $\Delta$ tensored with identity maps on its edges. After suitable degree shifts, by collapsing the cube (similar to passing to the total complex of a polycomplex) one obtains a complex $C(D)$ of graded abelian groups with a differential that preserves the quantum degree. Reidemeister moves can be lifted to specific homotopy equivalences between the complexes. Consequently, the isomorphism class of the bigraded homology groups $\Kh(L)\coloneqq\Kh(C(D))$ is an invariant of $L$, now widely called \emph{$\mathfrak{sl}_2$ homology} or \emph{Khovanov homology}. Identification of the Jones polynomial as the Euler characteristic of $\Kh(L)$ is immediate, since the construction of $C(D)$  lifts Kauffman's inductive formula. 

One can think of this construction of a  link homology as coming from a commutative Frobenius algebra $A$ over $\Z$, as above. The key property  of $A$ is having rank two over the ground ring $\Z$: using an algebra $A$ of larger rank, the homology fails to be invariant under Reidemeister I moves. On the other hand, a modification of this construction, deforming the relation $X^2=0$, gives rise to so-called \emph{equivariant} link homology~\cite{BN05:TangleCob,Khovanov06:Frobenius}. The essentially most general deformation comes from working over the ground ring $R'=\Z[h,t]$ and setting $A'$ to be
\begin{equation}\label{eq_frob_2}
  A' = R'[X]/(X^2-hX-t), \qquad \epsilon\co A' \to R', \quad  \epsilon(1)=0, \quad \epsilon(X)=1. 
\end{equation}
The equivariant theory turns out to be important for applications (see Section~\ref{sec:ssec} and~\ref{sec:apps}).

As mentioned above, this construction of link homology can be phrased via a rank two commutative Frobenius pair $(R,A)$, giving rise to a 2-dimensional topological quantum field theory (TQFT) $F=F_A$ with $F(\emptyset)=R$ (with is $\Z$ or $R'$ above) and $F(S^1)\cong A$. 
That a 2D TQFT of rank two can be bootstrapped into a link homology theory was surprising.  

There is also a \emph{reduced} version of the invariant, corresponding to the normalization $J(U)=1$. Fix a marked point on a strand of $D$. There is a subcomplex $\widetilde{C}(D)\subset C(D)$ where the marked circle is labeled $X$ throughout. Shifting the quantum grading of $\widetilde{C}(D)$ down by $1$ and taking homology gives $\rKh(L)$, the \emph{reduced Khovanov homology}. It is easy to see that $C(D)/\widetilde{C}(D)\cong \widetilde{C}(D)$, so there is a long exact sequence
\begin{equation}\label{eq:rKh-seq}
  \cdots\to \rKh^{i,j-1}(L)\to \Kh^{i,j}(L)\to \rKh^{i,j+1}(L)\to \rKh^{i+1,j-1}(L)\to\cdots.
\end{equation}

A paper of D.~Bar-Natan~\cite{BN02:expository} helped to provoke early interest in the subject, as well as giving computations of $\Kh(K)$ for knots through 12 crossings. (More work on computing $\Kh(L)$ is described in Section~\ref{sec:apps}.)

\subsection{Tangles and representations}
The Kauffman bracket invariant admits a relative version for tangles in the 3-disk~\cite{Kauffman87:state,Kauffman88:statistical,KauffmanLins94,CFS95:6j}. Start with a tangle $T$ in  $\mathbb{D}^3$ with $2n$ boundary points, and consider a generic projection of it to the 2-disk $\mathbb{D}^2$, with $2n$ boundary points spread out around the boundary $\partial\mathbb{D}^2$. Let $\Kaum_n$ be the free $\Z[q,q^{-1}]$-module with basis $B^n$ the set of $\emph{crossingless matchings}$ of $2n$ boundary  points via $n$ disjoint arcs inside a disk. The relative Kauffman bracket associates to $T$ an element $\brak{T}$ of $\Kaum_n$ by resolving each crossing following Kauffman's recipe. The braid group on $2n$ strands acts on $\Kaum_n$ by attaching a braid to a crossingless matching and then reducing the result via Kauffman's relations. (In fact, the larger group of annular braids acts). More generally, a tangle $T$ in a strip $\R\times [0,1]$ with $2n$ bottom and $2m$ top points (a $(2m,2n)$-tangle) induces a $\Z[q,q^{-1}]$-linear map 
\begin{equation}
    K(T) \co \Kaum_n \to \Kaum_m. 
\end{equation}
These maps fit together into a functor from the category of \emph{even} tangles (tangles with an even number of top and bottom endpoints) to the category of $\Z[q,q^{-1}]$-modules. Variations of Kauffman's construction can be made into monoidal functors from the category of tangles that assign $n\th$ tensor power of the fundamental representation $V$ of quantum $\mathfrak{sl}_2$ (or a suitable subspace of $V^{\otimes n}$) to $n$ points on the plane and intertwiners between tensor powers of representations to tangles. The above setup with crossingless matchings correspond to assigning the subspace of invariants $\mathrm{Inv}_{U_q(\mathfrak{sl}_2)}(V^{\otimes n})$ to $n$ points. This subspace is trivial
when $n$ is odd and has a basis of crossingless matchings for even $n$~\cite{KauffmanLins94,CFS95:6j,Khovanov97:thesis,FrenkelKhovanov97:bases}.  

\vspace{0.1in} 

\begin{figure}
  \centering
  \begin{tabular}{cccc}
    \quad
  \begin{tikzpicture}[scale=.8]
        \draw[dashed] (0,0) circle (1);
        \draw[thick] (.707,-.707) -- (-.707,.707);
        \draw[thick] (-.707,-.707) -- (-.1,-.1);
        \draw[thick] (.1,.1) -- (.707,.707);
      \end{tikzpicture} \quad & \quad 
\begin{tikzpicture}[scale=.8]
  \draw[dashed] (0,0) circle (1);
  \draw[-, thick] (-.707,-.707) .. controls (0,-.1) .. (.707,-.707);
  \draw[-, thick] (-.707,.707) .. controls (0,.1) .. (.707,.707);
\end{tikzpicture}   \quad & \quad
\begin{tikzpicture}[scale=.8]
    \draw[dashed] (0,0) circle (1);
    \draw[-, thick] (-.707,-.707) .. controls (-.1,0) .. (-.707,.707);
    \draw[-, thick] (.707,-.707) .. controls (.1,0) .. (.707,.707);
  \end{tikzpicture}
  \quad  & \quad \includegraphics[width=.75in]{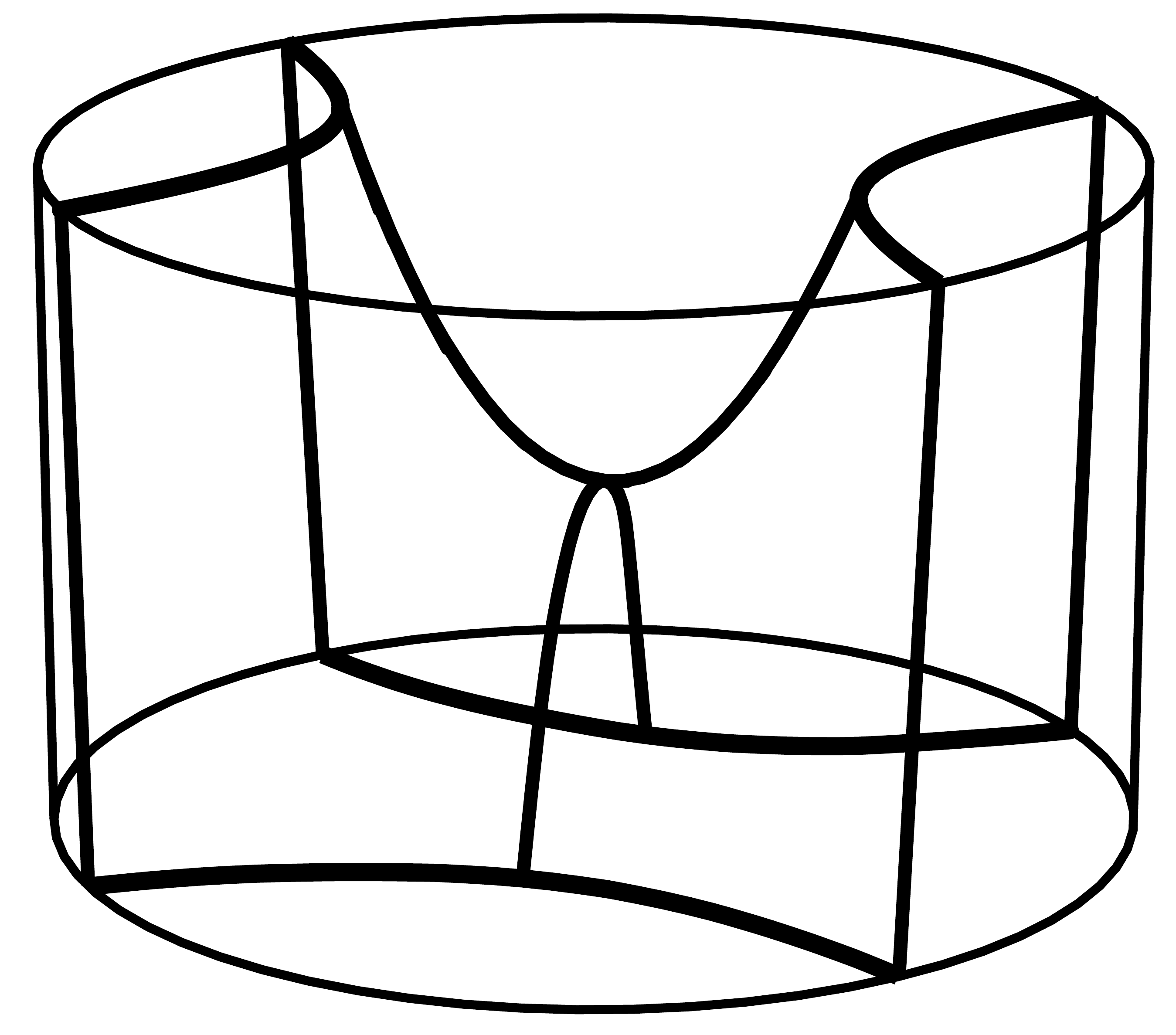}\quad{}
    \\[.5em]
    $T$ & $a$ & $b$ & $s$
    \end{tabular}
  \caption{\textbf{The complex associated to a tangle.} The complex of graded projective $H^2$-modules $F(T)$ is given by
    $ 0\to P_a\{1\}\stackrel{F(s)}{\longrightarrow} P_b\to 0
      $
    where $P_a$ and $P_b$ are the modules associated to the flat
    tangles $a$ and $b$ shown, and $s$ is the indicated
    saddle cobordism. The notation $\{1\}$ indicates a quantum grading shift.
  }
  \label{fig:tangle-cx}
\end{figure}
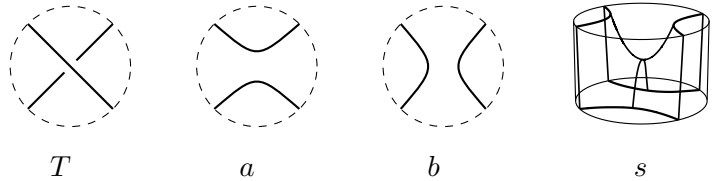

Upon categorification, $\Kaum_n$ becomes a Grothendieck group of a suitable category $\mcC_n$. A crossingless matching $a\in  B^n$ with $2n$ specified endpoints $\underline{p}$ becomes an object $P_a$ of $\mcC_n$. We can guess that hom spaces $\Hom_{\mcC_n}(P_a,P_b)$ will come from cobordisms between $a$ and   $b$, that is, surfaces $S$ embedded in $\mathbb{D}^2\times [0,1]$ with boundaries $a\times\{0\}$, $b\times\{1\}$, and $[0,1]\times\underline{p}$. (An example is on the right of Figure~\ref{fig:tangle-cx}.) The  total boundary  of  such surface $S$ is homeomorphic to the  $1$-manifold $\overline{b}a$ given by gluing $a$ and $b$ along their boundary points. One can then define 

\begin{equation*}
\Hom_{\mcC_n}(P_a,P_b) \coloneqq \ F(\overline{b}a) 
\end{equation*}
by applying the 2D TQFT $F$ as  above to that $1$-manifold. It is straightforward to define associative multiplications
\begin{equation*}
\Hom_{\mcC_n}(P_a,P_b) \times   \Hom_{\mcC_n}(P_b,P_c) \ \lra \ 
\Hom_{\mcC_n}(P_a,P_c), \ \ a,b,c\in B^n
\end{equation*}
by applying $F$ to appropriate cobordisms~\cite{Kh02:tangles}.

More carefully, to define $\mcC_n$ we start with  
objects $\{P_a\}_{a\in B^n}$ and morphisms as above and form a pre-additive category $\mcC_n''$. Equivalently, category $\mcC_n''$ can be viewed as an idempotented ring 
\begin{equation*}
    H^n \coloneqq \oplusop{a,b\in B^n} \Hom_{\mcC_n''}(P_a,P_b) = \oplusop{a,b\in B^n} F(\overline{b}a),
  \end{equation*}
  the \emph{arc ring},
with idempotents $1_a\in F(\overline{a}a)$ given by identity cobordisms from $a$ to itself. It is also possible to keep track of morphisms in different degrees and refine the category by restricting morphisms to degree $0$ parts of graded abelian groups $F(\overline{b}a)$ but allowing  grading shifts of generating objects to capture the entire groups. 

From the idempotents $1_a$ one can recover the projective modules $P_a\coloneqq H^n\, 1_a$  over $H^n$.

One can then form an additive closure $\mcC_n'$ of  the category  $\mcC_n''$ by also allowing finite direct sums of objects. The category $\mcC_n'$ happens to be Karoubi closed, which is not  hard to check and simplifies working with it. The category $\mcC_n'$ is equivalent to the category of graded projective finitely-generated modules over the graded ring $H^n$. 

To a flat (crossingless) tangle $T$ in a disk $\mathbb{D}^2$ with $2n$ endpoints there is associated an  object $F(T)$ of $\mcC_n'$ or, equivalently, a projective graded $H^n$-module. If $T$ is the union of $k$ circles and a crossingless matching $a\in B^n$, then the projective module is isomorphic  to $A^{\otimes k}\otimes P_a$, that  is, to the sum of  $2^k$ copies of the projective module $P_a$, with appropriate grading shifts.

The Grothendieck group $K_0(\mcC'_n)$ of $\mcC'_n$ is a free $\Z[q,q^{-1}]$-module with basis given by the symbols $[P_a]$ of projective modules, over all crossingless matchings $a\in B^n$. This Grothendieck group can also be defined as $K_0$ of the graded algebra $H^n$. There is a canonical isomorphism of $\Z[q,q^{-1}]$-modules 
\begin{equation*}
     K_0(\mcC'_n) \ \cong \ \Kaum_n.
\end{equation*}
Now form the category $\mcC_n$ of bounded  complexes of objects of $\mcC_n'$, modulo chain homotopies. The inclusion $\mcC_n'\subset \mcC_n$ induces an isomorphism of their Grothendieck groups. 

To a planar diagram $D$ of a tangle $T$ with $2n$ endpoints there is an associated object object $F(D)$ of $\mcC_n$, by a relative version of the cube construction. 
Namely, define $F(D)$ to be the iterated mapping cone of the two resolutions at each crossing, that is, the total complex of the cube of resolutions of $D$. See Figure~\ref{fig:tangle-cx} for a simple example.

Reidemeister moves of tangle diagrams lift to chain homotopy equivalences, and the isomorphism class of the object $F(D)$ is an invariant of $T$. On the Grothendieck group, $F(D)$ descends to the element $\brak{T}\in\Kaum_n$. 

Similarly, given a tangle diagram $D$ with $2m$ bottom and $2n$ top endpoints  there is an associated complex of $(H^m,H^n$)-bimodules, and tensoring with this complex of bimodules gives an exact functor $F(D)\co\mcC_n\to \mcC_m$. This construction lifts to a 2-functor from the category of flat tangles and their cobordisms to the category of bimodules and their homomorphisms.  Furthermore, it lifts to a projective functor (well-defined on 2-morphisms up to an overall sign) from the 2-category of tangle cobordisms to the 2-category of complexes of bimodules over $H^n,$ over all $n\ge 0$, and maps of complexes, up to homotopy~\cite{Kh06:cobordisms,BN05:TangleCob} (see also~\cite{Jacobsson} for another proof). Taking care of the sign is subtle; see~\cite{CMW09:signs,Caprau08:signs,Blanchet10:signs,Sano:signs}.

\vspace{0.1in} 

Categories of representations of the arc rings $H^n$ categorify $\Kaum_n\cong \mathrm{Inv}_{U_q(\mathfrak{sl}_2)}(V^{\otimes n})$. It turns out that the entire tensor product $V^{\otimes n}$, as well as the commuting actions of the Temperley-Lieb algebra and quantum $\mathfrak{sl}_2$ on it, can also be categorified. This categorification  was realized in~\cite{BFK} via maximal singular and parabolic blocks of highest weight categories for $\mathfrak{sl}_n$, with the commuting actions lifting to those by projective functors and Bernstein-Zuckerman functors (see also~\cite{FKS06,Zheng08}). 

\vspace{0.1in} 

The tensor power $V^{\otimes n}$ decomposes as the sum of its weight spaces $V^{\otimes n}(k)$, $k=0,\dots, n$. 
A more explicit categorification of weight spaces and the Temperley-Lieb algebra action on them can be achieved via specific subquotient rings of $H^{2n}$~\cite{BS11:I,ChenKhovanov14}. J.~Brundan and C.~Stroppel  showed~\cite{BS11:III,BS12:IV} that these subquotient rings (a) describe maximal  parabolic blocks of highest weight categories for $\mathfrak{sl}_n$, relating the two categorifications, and (b) describe blocks of representations of Lie superalgebras $\mathfrak{gl}(m|k)$.

The space $\mathrm{Inv}_{U_q(\mathfrak{sl}_2)}(V^{\otimes n})$ of invariants is naturally a subspace of the middle weight space $V^{\otimes 2n}(n)$. Analogues of this subspace for a general weight space $V^{\otimes n}(k)$ are given by the kernel of the generator $E\in \mathfrak{sl}_2$ for $k\ge n/2$, and the kernel of $F\in \mathfrak{sl}_2$ for $k\le n/2$. Categorifications of these subspaces are provided by representation categories of certain Frobenius algebras, like $H^n$, that can be obtained as subquotients of $H^n$. The latter Frobenius algebras and  Morita and derived Morita equivalent algebras are widespread in modular representation theory.  For instance, Hiss and Lux's book~\cite{HissLux89:book} lists hundreds of examples of blocks of finite groups over finite characteristic fields that are (derived) Morita equivalent to the self-dual part of the zigzag algebra from~\cite{KhovanovSeidel02}, the latter giving a categorification of the reduced Burau representation of the braid group and of the corresponding subspace of the first nontrivial weight space, $V^{\otimes n}(1)$. 

\vspace{0.1in} 

A very general framework for a categorification of tensor products of quantum group representations and Reshetikhin-Turaev link invariants was developed by Ben Webster~\cite{Webster17:book}. The $\mathfrak{sl}_2$ case of his construction~\cite{Webster16:Grass} uses algebras that are Morita equivalent to Koszul duals of the abovementioned subquotients of $H^{2n}$.

\subsection{Connections to algebraic geometry, symplectic geometry,
  and beyond}
The connection with representation theory inspired a further connection with symplectic geometry. Given a symplectic manifold $(M,\omega)$, there is an associated triangulated category, the derived Fukaya category; the objects of the Fukaya category are Lagrangian submanifolds of $M$ (with certain extra data) and the morphism spaces are categorified intersection numbers, defined via Floer theory. Given a braid group action on $(M,\omega)$, there is an induced braid group action on the Fukaya category and hence, potentially, a knot invariant. The first examples of such braid group actions were given by P. Seidel and the first author~\cite{KhovanovSeidel02}. Soon after, Seidel and I.~Smith gave a more braid group action on a more complicated, but natural, symplectic manifold, and from it a conjectured Floer-theoretic definition of Khovanov homology, which they called \emph{symplectic Khovanov homology}~\cite{SS06:Kh-symp}. (See also~\cite{Manolescu06:Hilbert} for a reinterpretation of this construction.) Recently, M.~Abouzaid and Smith proved that this conjecture holds over $\Q$~\cite{AS16:formal,AS19:KhFloer}. The proof uses the extension of Khovanov homology to tangles discussed above to identify the two theories. At present, it is unknown whether the torsion in symplectic Khovanov homology and in combinatorial Khovanov homology agree. Although it is harder to compute, symplectic Khovanov homology is in some ways more geometric. In particular, its relationship to Heegaard Floer homology and its behavior for periodic knots (see Section~\ref{sec:ssec}), as well as the equivariant versions of the theory in the sense of~\eqref{eq_frob_2} above, all have geometric definitions via group actions on the symplectic manifold~\cite{SS10:localization,HLS16:equi}.

The symplectic manifolds in the  Seidel-Smith construction are examples of quiver varieties, so carry hyperk\"ahler structures. Complex Lagrangians determine objects of both the Fukaya category and
the category of coherent sheaves with respect to the rotated almost
complex structure. The fact that the automorphism algebras on the two
sides are isomorphic to ordinary cohomology can be seen as a shadow of
mirror symmetry, and can often be lifted to an equivalence of
categories. Consequently, one would expect that the tangle extension of Khovanov homology can be realized via derived categories of coherent sheaves on the corresponding quiver varieties, with functors associated to tangles acting via suitable Fourier-Mukai kernels (convolutions with objects of the derived category on the direct product of varieties). A modification of this idea was realized by S.~Cautis and J.~Kamnitzer~\cite{CautisKamnitzer08:sl2}. They use certain smooth completions of these quiver varieties which can be realized  as iterated $\mathbb{P}^1$-bundles and interpreted as convolution varieties of the affine Grassmannian for $SL(2)$, also providing a connection to the geometric Satake correspondence. The relation to quiver varieties and the $(n,n)$-Springer fiber has been established by R.~Anno~\cite{Anno:affine} and Anno and V.~Nandakumar~\cite{AnnoNandakumar}, who also explained the relation between coherent sheaves on these varieties and the rings $H^n$ and their annular versions.
An isomorphism between the center of $H^n$ and the cohomology ring of the $(n,n)$-Springer fiber, established in~\cite{Kh04:Springer}, was an earlier indication of the connection between the two structures. 

\vspace{0.1in} 

There has been strong interest in giving physical reinterpretations and extensions of link homology invariants. One program to do so was initiated by Witten, using the Kapustin-Witten and Haydys-Witten equations~\cite{Witten12:Kh}. Other proposals have been put forward by S.~Gukov, A.~Schwarz, and C.~Vafa~\cite{GSV05:KR}, Gukov, P.~Putrov, and Vafa~\cite{GPV17:fivebranes}, Gukov-D. Pei-Putrov-Vafa~\cite{GPPV20:BPS}, M.~Aganagic~\cite{Aganagic:1,Aganagic:2}, and others.

Currently, Khovanov homology is only defined for links in a few manifolds: $S^3$, as described above; links in thickened surfaces, in work of Asaeda-Przytycki-Sikora~\cite{APS04:annular}; and links in connected sums of $S^2\times S^1$, in work of Rozansky~\cite{Rozansky:S2S1} and Willis~\cite{Willis21:S2S1}. (See also the universal construction in~\cite{MWW:blob}.) One appeal of some of the conjectural physical approaches to Khovanov homology is that they may apply in general 3-manifolds. 
In a recent paper~\cite{QiSussan:pdg},  J.~Sussan and Y.~Qi  categorify the Jones polynomial when the quantum parameter $q$ is a prime root of unity; this is also related to extending Khovanov homology to other 3-manifolds. 

There is a large literature on categorification of $\mathfrak{sl}(k)$ representations and quantum invariants,  for an arbitrary $k$. For lack of space, we will not discuss these developments in this paper. Nor do we discuss the related topics of annular homology, categorifications of the colored Jones polynomial, foams, and categorified quantum groups.

\section{Signs and spectral sequences}\label{sec:ssec}
One reason Khovanov homology has been important is that it seems to be a kind of free object in the category of knot homologies, a property which is witnessed by the many spectral sequences from Khovanov homology to other knot homologies. (An attempt to make precise the sense in which Khovanov homology is free was given in~\cite{BHL19:Kh-Floer}.) These spectral sequences often connect invariants whose constructions appear quite different, in some cases giving relationships between invariants that are not apparent at the classical, decategorified level. They have led to many of the topological applications of Khovanov homology, as well as to new properties of Khovanov homology itself.

The first spectral sequence from Khovanov homology was constructed by E.~S.~Lee~\cite{Lee05:alternating} (see also~\cite{Rasmussen10:s}). Recall the family of deformations $A'$ of the Frobenius algebra $A$ from Equation~\eqref{eq_frob_2}. Taking the parameters $(h,t)=(0,1)$ and extending scalars from $\Z$ to $\Q$, we obtain the algebra $\Q[X]/(X^2=1)$. The quantum grading weakens to a filtration on the resulting complex, inducing a spectral sequence from Khovanov homology to this deformed knot invariant, called \emph{Lee homology}. To understand Lee homology, note that this Frobenius algebra diagonalizes, as a direct sum of two one-dimensional Frobenius algebras. It follows easily that the Lee homology of a $c$-component link has dimension $2^{c}$. Using this construction, Lee verified a conjecture of Bar-Natan~\cite{BN02:expository}, S.~Garoufalidis~\cite{Gar04:conjecture}, and the first author~\cite{Khovanov03:patterns} that the Khovanov homology of an alternating knot lies on two adjacent diagonals. More famous applications of this spectral sequence are discussed in Section~\ref{sec:apps}.

Most of the other spectral sequences from Khovanov homology relate to gauge theory. The first of these is due to P.~Ozsv\'ath and Z.~Szab\'o~\cite{OS05:dcov}. Given a closed, oriented $3$-manifold $Y$, they had constructed an abelian group $\HFa(Y)$, the homology of a chain complex $\CFa(Y)$~\cite{OS04:HolomorphicDisks}. Inspired  by A.~Floer's exact triangle~\cite{Floer95:triangle}, they showed that given a knot $L\subset Y$ and slopes $\mu$, $\lambda$, and $\mu+\lambda$ on $\partial(\mathrm{nbd}(L))$ intersecting each other pairwise once, there is an exact triangle relating the Floer homologies of the surgeries $\HFa(Y_\mu(L))$, $\HFa(Y_\lambda(L))$, and $\HFa(Y_{\mu+\lambda}(L))$~\cite{OS04:HolDiskProperties}. In particular, given a link $K$ in $S^3$, if $K_0$ and $K_1$ are the $0$- and $1$-resolutions of a crossing of $K$ then the surgery exact triangle gives an exact triangle of Floer homologies of their branched double covers,
\[
  \begin{tikzpicture}
    \node at (0,0) (L0) {$\HFa(\Sigma(K_0))$};
    \node at (3,0) (L1) {$\HFa(\Sigma(K_1))$};
    \node at (1.5,-1) (L) {$\HFa(\Sigma(K))$};
    \draw[->] (L0) to (L1);
    \draw[->] (L1) to (L);
    \draw[->] (L) to (L0);
  \end{tikzpicture}
\]
(This is an ungraded exact triangle: the groups $\HFa(\Sigma(K))$ do not have canonical $\Z$-gradings, and the gradings they do have are not respected by the maps in the exact triangle.) The surgery exact triangle is local, in the sense that given disjoint links $L$ and $L'$, the maps in the surgery exact triangles associated to $L$ and $L'$ commute or, at the chain level, commute up to reasonably canonical homotopy. So, resolving all $N$ crossings of $K$ gives a cube of resolutions for $\CFa(\Sigma(K);\FF_2)$. The $E^1$-page of the associated spectral sequence is
\begin{align*}
  \bigoplus_{I\in 2^{N}}\HFa(\Sigma(K_I);\FF_2)&=\bigoplus_{I\in 2^N}\HFa(\#^{|K_I|-1}(S^2\times S^1);\FF_2)\\
  &=\bigoplus_{I\in 2^N}\HFa(S^2\times S^1;\FF_2)^{\otimes(|K_I|-1)}\cong (\FF_2\oplus \FF_2)^{\otimes|K_I|-1},
\end{align*}
which has the same dimension as the reduced Khovanov complex. The differential on the $E^1$-page comes from merge and split cobordisms $(S^2\times S^1) \leftrightarrow (S^2\times S^1)\#(S^2\times S^1)$. These maps correspond to some 2-dimensional Frobenius algebra which, in fact, turns out to be the algebra $A$. Thus, one obtains a spectral sequence from the reduced Khovanov homology of (the mirror of) $K$, with $\FF_2$-coefficients, to $\HFa(\Sigma(K);\FF_2)$.

The Euler characteristic of $\HFa(\Sigma(K))$ is the number of elements in $H_1(\Sigma(K))$ if finite, or $0$ otherwise. So, the Ozsv\'ath-Szab\'o spectral sequence lifts the equality  $J_{-1}(K)=\det(K)$. 

To summarize, the key properties of $\HFa(\Sigma(K))$ used to construct the Ozsv\'ath-Szab\'o spectral sequence were the existence of an unoriented skein triangle satisfying a far-commutativity property, TQFT properties for disjoint unions, merges, and splits, and the fact that its value on an unknot (or, more accurately, 2-component unlink) is a 2-dimensional vector space.

In 2010, P.~Kronheimer and T.~Mrowka built a gauge-theoretic invariant $I^\natural$ with these properties, using Donaldson theory~\cite{KM11:unknot}. Like many gauge-theoretic invariants, the value of $I^\natural$ constrains how surfaces can be embedded. Using this, Kronheimer-Mrowka deduced that if the genus of a knot $K$ is $>1$ then $I^\natural(S^3,K)$ has dimension $>1$. From the argument above, there is a spectral sequence $\rKh(K)\Rightarrow I^\natural(S^3,K)$, hence:
\begin{theorem}\label{thm:KM-detect}(Kronheimer-Mrowka \cite{KM11:unknot})
  If $\rank\, \rKh(K)=1$ then $K$ is the unknot.
\end{theorem}
The stronger, and older, conjecture, that $J(K)=q+q^{-1}$ only if $K$ is the unknot, remains open.

There are many other spectral sequences from Khovanov homology, including more variants of the Lee spectral sequence~\cite{DGR06:super,Ballinger}, spectral sequences defined using instanton and monopole Floer homology~\cite{Daemi,Scaduto15:ss,Bloom11:ss}, other spectral sequences defined via variants of Heegaard Floer homology~\cite{GW10:ss,Roberts13:ss}, 
spectral sequences coming from equivariant symplectic Khovanov homology and equivariant Khovanov homology~\cite{SS10:localization,Cornish,Zhang18:periodic,SZ:periodic}, and a combinatorial spectral sequence conjectured to agree with the Ozsv\'ath-Szab\'o spectral sequence~\cite{Szabo15:geometric} (see also~\cite{SSS17:perturb}). This last spectral sequence also supports another conjecture: that the Ozsv\'ath-Szab\'o spectral sequence preserves the \emph{$\delta$-grading} $\delta=j-2i$ on Khovanov homology~\cite{Greene13:spanning-tree}. Another notable spectral sequence is due to J.~Batson and C.~Seed: given a link $L=L_1\cup L_2$, they construct a spectral sequence $\Kh(L_1\cup L_2)\Rightarrow \Kh(L_1\amalg L_2)$  to the disjoint union of the sub-links $L_1$ and $L_2$~\cite{BS15:splitting} (which is just $\Kh(L_1)\otimes \Kh(L_2)$ if working over a field). The page of collapse of this spectral sequence  gives a lower bound on the unlinking number of $L$. It and many of the other spectral sequences have also been used to prove further detection results for Khovanov homology, in the spirit of Theorem~\ref{thm:KM-detect}. Often, the proofs of detection results combine several of these spectral sequences. Some examples of such results include:

\begin{theorem}(Batson-Seed~\cite{BS15:splitting})
  Let $U^m$ be the $m$-component unlink. If $\dim \Kh^{i,j}(L;\FF_2)= \dim\Kh^{i,j}(U^m;\FF_2)$ for all $i$ and $j$ then $L$ is isotopic to $U^m$.
\end{theorem}
The proof uses Theorem~\ref{thm:KM-detect} and the Batson-Seed
spectral sequence. A related result was obtained earlier by M.~Hedden
and Y.~Ni~\cite{HN13:detection}. (By contrast, the Jones polynomial does not detect the unlink~\cite{Thistlethwaite01:JonesUnlink,EKT03:JonesUnlink}. Indeed, most of the detection results mentioned below also do not hold for the Jones polynomial.)

\begin{theorem}(Xie-Zhang \cite{XZ:link-detect})
  If $K$ is an $m$-component link with $\dim\Kh(K;\FF_2)=2^m$ then $K$ is a forest of Hopf links.
\end{theorem}
The proof uses Kronheimer-Mrowka's spectral sequence and its extension to annular links \cite{Xie:annular} (building on~\cite{APS04:annular,Roberts13:annular,GW10:ss}); Batson-Seed's spectral sequence; and N.~Dowlin's spectral sequence mentioned below. In other papers, the authors classify all links with Khovanov homology of dimension $\leq 8$~\cite{XZ:small} and show that Khovanov homology detects, for instance, $L7n1$~\cite{XZ:L7n1}. Similarly, 
Khovanov homology detects the link $T(2,6)$~\cite{Martin:T26}.

\begin{theorem} Let $K$ be a knot.
  \begin{enumerate}
  \item (Baldwin-Sivek \cite{BS:trefoils}) If $\Kh(K)\cong \Z^4\oplus\Z/2\Z$ then $K$ is the trefoil knot. 
  \item (Baldwin-Dowlin-Levine-Lidman \cite{BDLLS}) If $\rank(\rKh(K))=5$ and the reduced Khovanov homology is supported in $\delta$-grading $0$ then $K$ is the figure 8 knot.
  \item (Baldwin-Hu-Sivek \cite{BHS:cinquefoil}) If $\Kh(K)\cong\Kh(T(2,5))$ then $K$ is the torus knot $T(2,5)$.
  \end{enumerate}
\end{theorem}
The proof of the first statement uses Kronheimer-Mrowka's spectral sequence, the second uses Dowlin's spectral sequence, and the third uses an annular version of the Kronheimer-Mrowka spectral sequence~\cite{Xie:annular,XZ:sut-tan}, Dowlin's spectral sequence, the spectral sequences for periodic knots~\cite{SZ:periodic,BPS:periodic} mentioned above, further hard results on Floer homology~\cite{Xie:annular,XZ:sut-tan,KLT20:HFHM,LeeTaubes12:periodic,CottonClay}, and the $\mathfrak{sl}_2(\C)$-action on annular Khovanov homology~\cite{GLW18:sl2}.

\vspace{0.1in} 

Some of these, like the spectral sequences from equivariant Khovanov homology, lift, or at least recall, well-known properties of the Jones polynomial, such as K.~Murasugi's formula~\cite{Murasugi88:JonesPeriodic}. By contrast, other spectral sequences seem invisible to the Jones polynomial. Perhaps most striking, building on work of Ozsv\'ath-Szab\'o, Ozsv\'ath-A.~Stipsicz-Szab\'o, and C.~Manolescu~\cite{OSz09:Cube,OSSz09:singular,Manolescu14:cube}, Dowlin showed~\cite{Dowlin:ss} that there is a spectral sequence from Khovanov homology to Heegaard Floer knot homology (which categorifies the Alexander polynomial). This implies that the dimension of Khovanov homology is always at least as large as that of knot Floer homology, a statement with no known analogue in terms of the classical Jones and Alexander polynomials (though see~\cite{GuManion}).

The reader might notice the prevalence of $\FF_2$-coefficients in these spectral sequences. In fact, many of the spectral sequences have lifts to $\Z$-coefficients, but do not start from Khovanov homology. Instead, they start from a variant, \emph{odd Khovanov homology}, discovered  by Ozsv\'ath-J.~Rasmussen-Szab\'o when trying to lift the Ozsv\'ath-Szab\'o spectral sequence to $\Z$-coefficients~\cite{OSSz13:odd-kh}. In constructing the cube of resolutions for Khovanov homology, to a collection of $n$ circles $Z_1,\dots,Z_n$ in the plane one associates $\bigl(\Z[X]/(X^2)\bigr)^{\otimes n}$, a quotient of the symmetric algebra on $n$ variables. To construct odd Khovanov homology, one instead associates the exterior algebra on $n$ variables, $\Lambda\langle Z_1,\dots,Z_n\rangle$. Merging circles $Z_i$ and $Z_j$ into a circle $Z$ corresponds to the map sending $Z_i$ and $Z_j$ to $Z$, while splitting $Z$ into $Z_i$ and $Z_j$ corresponds to  multiplying by $(Z_i-Z_j)$ (or $(Z_j-Z_i)$: the definition involves a choice, which one can fix by picking certain orientations at the crossings). The resulting cube neither commutes nor anti-commutes, but nonetheless one can show that it is possible to assign signs to the edges, in an essentially unique way, to get an anti-commuting cube. The homology of the total complex of this cube is odd Khovanov homology.

Although the change to the definition might seem slight, odd Khovanov homology has quite different properties from ordinary Khovanov homology:
\begin{itemize}
\item Unreduced odd Khovanov homology is the direct sum of two copies of reduced odd Khovanov homology, while for ordinary, even Khovanov homology the long exact sequence~\eqref{eq:rKh-seq}
  almost never splits.
\item Odd Khovanov homology is mutation invariant~\cite{Bloom10:mutation}, while even Khovanov homology of links is not~\cite{Wehrli:not-mutation}.
\item There is no known analogue of the Lee spectral sequence for odd Khovanov homology, but rather there is a spectral sequence from reduced odd Khovanov homology to $\Z^{\det(L)}$~\cite{Daemi}.
\end{itemize}
For alternating knots, by the first point above, odd Khovanov homology has no torsion, while even Khovanov homology almost always has $2$-torsion, but no other torsion~\cite{Shu14:torsion,Shu:2torsion}. On the other hand, more torsion appears in the reduced odd Khovanov homology than in reduced even Khovanov homology for small knots~\cite{Shu11:odd}. (See also~\cite{LowSaz17:chromatic,MPSWYY18:torsion} for further results and citations.)

The representation-theoretic interpretation of odd Khovanov homology is substantially more involved than that of Khovanov homology, and is still an active area of research (see, for instance,~\cite{Putyra14:chronology,NP:odd,EllisQi16:odd,LaudaRussel14:odd}).

\section{Spectrification}\label{sec:spectra}
As we saw above, Khovanov homology is closely related to
low-dimensional Floer homologies, a family of invariants defined using
a kind of semi-infinite-dimensional Morse theory. Unlike
R.~Palais and S.~Smale's infinite-dimensional Morse theory~\cite{PalaisSmale64},
Floer homology is not isomorphic to the singular homology of the
ambient space.  R.~Cohen, J.~Jones, and G.~Segal proposed an alternate construction
of a stable homotopy type, or spectrum, $X$ associated to a Floer homology setup so
that the (reduced) homology of $X$ is isomorphic to the Floer
homology under consideration. (Unlike ordinary cohomology, Floer cohomology rarely has a
graded-commutative cup product, so it is natural to expect it would be
associated to a spectrum rather than a space.) Cohen-Jones-Segal's
original construction has only been made rigorous in a few cases but,
using other techniques, Manolescu did construct a stable homotopy
refinement of Seiberg-Witten Floer
homology~\cite{Manolescu03:SW-htpy}. Given Seidel-Smith's conjectured
Floer homology formulation of Khovanov homology~\cite{SS06:Kh-symp},
Khovanov homology's close relationship to Seiberg-Witten Floer
homology~\cite{OS05:dcov,Bloom11:ss}, and Manolescu's stable homotopy
refinement of Seiberg-Witten Floer homology, it was natural to expect
that there would be a stable homotopy refinement of Khovanov homology,
and in fact S.~Sarkar and the second author showed that there is~\cite{LS14:kh-htpy}. Another construction of such
a stable homotopy type was soon given by
P.~Hu, D.~Kriz, and I.~Kriz~\cite{HKK16:Kh-htpy}; somewhat later, the two
constructions were shown to be
equivalent~\cite{LLS20:khovanov-product}.

The idea behind Cohen-Jones-Segal's construction is as follows. First,
consider building a CW complex from a $0$-cell, an $n$-cell, and an
$(n+k)$-cell. The attaching data for the $(n+k)$-cell is a map $S^{n+k-1}\to S^n$. If
$n>k+1$ then, by the Pontrjagin-Thom construction, this is equivalent
to specifying a manifold $M^{k-1}$ and a framing of its stable normal
bundle. Next, suppose we want to build a space from cells of dimension
$0$, $n$, $n+k$, and $n+k+\ell$, where $n$ is large compared to $k$
and $\ell$. One can specify the attaching map $\partial e^{n+k}\to S^n$ by
a closed, framed manifold $M^{k-1}$, and the quotient of the attaching
map $\partial e^{n+k+\ell-1}\to X^{n+k}/X^{n+k-1}=e^{n+k}/\partial e^{n+k}=S^{n+k}$
by a closed, framed manifold $N^{\ell-1}$. It is not hard to see that
this map factors through a map $S^{n+k+\ell-1}\to X^{n+k}$ if and only
if the product $M\times N$ is the boundary of a framed manifold
$P^{k+\ell-1}$, and a choice of such a lift  up to homotopy is the same
as a choice of $P^{k+\ell-1}$ up to appropriate framed cobordism. Continuing this line of reasoning to an
arbitrary number of cells leads to the notion of a \emph{framed flow
  category}: a non-unital category where the morphism spaces are
manifolds with corners, and the composition maps sweep out the
boundaries. Such a framed flow category specifies a CW complex, its
\emph{realization}, by a Pontrjagin-Thom construction as above.

A Morse function $f$ on a closed manifold $M$, together with a generic
Riemannian metric, specifies a framed flow category with objects
corresponding to the critical points of $f$. The morphism space from
$x$ to $y$ is the space of gradient flow lines of $f$ from $x$ to
$y$. For example, for the usual Morse function on the circle the flow
category has two objects $S$ and $N$, and $\Hom(N,S)$ consists of two
points (with opposite framings). The flow category for the product
Morse function on the torus $S^1\times S^1$ has four critical points,
$SS$, $SN$, $NS$, and $NN$. The morphism sets $\Hom(SN,SS)$,
$\Hom(NS,SS)$, $\Hom(NN,SN)$, and $\Hom(NN,NS)$ consist of two points
each. The space $\Hom(NN,SS)$ is a disjoint union of four
intervals. For $S^1\times S^1\times S^1$, the product flow category
has $\Hom(NNN,SSS)$ a disjoint union of hexagons. For $S^1\times
S^1\times S^1\times S^1$, the product flow category has
$\Hom(NNNN,SSSS)$  a disjoint union of 3-dimensional permutohedra, and
now we have the general pattern: the flow category of $T^n$ is built
from permutohedra of dimension $0,\dots,n-1$. (The appearance of
permutohedra is not special to tori: they appear in any product
$X_1\times\cdots\times X_n$.)

With this in mind, specifying a stable homotopy refinement of Khovanov homology is
equivalent to building a framed flow category whose objects correspond
to the generators of the Khovanov complex. The morphism sets between
generators in adjacent gradings are framed $0$-manifolds, and counting
the number of points in these $0$-manifolds should give the
coefficients in the differential on the Khovanov complex. It turns out
not to be hard to define such a framed flow category, where all the
morphism sets are modeled on disjoint unions of permutohedra
(corresponding, perhaps, to all the tensor products appearing in the
Khovanov complex). So, like the Khovanov complex itself, in some sense
this appears to be the simplest, or freest, possible construction. In
fact, at present, it is equivalent to all known constructions of
functorial stable homotopy refinements of Khovanov homology.

Like Khovanov homology and the Jones polynomial itself, this
homotopy refinement of Khovanov homology is not built intrinsically
from a knot, but rather inductively via the cube of resolutions. So,
one must check that, up to stable homotopy equivalence, the result is 
independent of the knot diagram. This turns out to be easy, via
Whitehead's theorem: all one needs to do is construct maps of spectra
inducing the usual isomorphisms on Khovanov homology. Since the
invariance proof for Khovanov homology boils down to repeatedly taking
subcomplexes and quotient complexes, lifts to the stable homotopy type
come for free.

A stable homotopy refinement induces Steenrod operations on Khovanov
homology. If the Khovanov homology has a sufficiently simple form then
these, in turn, determine the stable homotopy type. The operation
$\mathrm{Sq}^1$ is just the Bockstein homomorphism, and one can give
an explicit formula for the Steenrod squaring operation
$\mathrm{Sq}^2$~\cite{LS14:square}, and more complicated formulas for
all Steenrod squares~\cite{Cantero20:squares}. The operation
$\mathrm{Sq}^2$ is enough to determine the stable homotopy type for
all knots up to 14 crossings and, in fact, some pairs of knots with
isomorphic Khovanov homologies are distinguished by their Steenrod
squares~\cite{Seed:square}. By introducing simplification operations
for flow categories, one can give computer computations for some more
complicated knots, and even by-hand computations for simple knots with nontrivial $\mathrm{Sq}^2$ operations, like
the $(3,4)$ torus knot~\cite{JLS17:moves,JLS20:calculus}.

Many structures for Khovanov homology can be lifted to or enhanced by
the stable homotopy type, including the Rasmussen
invariant~\cite{LS14:s,LLS20:khovanov-product} (see
Section~\ref{sec:apps}), Plamenevskaya's transverse
invariant~\cite{LNS15:transverse} (again, see
Section~\ref{sec:apps}), and the arc
algebras~\cite{LLS-kh-tangles} (see Section~\ref{sec:Kh}). There is an
analogue for odd Khovanov homology~\cite{SSS20:odd} (see
also~\cite{PS16:Bockstein}). The homotopical
refinement can even be used to prove new results about Khovanov homology
itself, such as formulas relating the Khovanov homology of periodic
links and their quotients~\cite{SZ:periodic,BPS:periodic}, partially
lifting results of Murasugi~\cite{Murasugi88:JonesPeriodic}.

While there has been some work on connections between these spectral
refinements and representation theory~\cite{AKW:sl2-spectrum,HKS19:derived}, even though the
original inspiration for the refinements comes from Floer theory,
direct connections with symplectic geometry or algebraic geometry remain unknown. There has also been work on giving stable homotopy refinements of $\mathfrak{sl}_n$ Khovanov-Rozansky homology, and connections between that and equivariant algebraic topology~\cite{Kitchloo,JLS19:matched}.

\begin{remark}
  These refinements are spectra $X$ whose singular homology is equal
  to Khovanov homology. The problem of finding a homotopy type $X$
  whose homotopy groups agree with Khovanov homology was also
  considered~\cite{ET14:htpy}. Unlike the case of homology, for
  homotopy groups there is a universal, functorial construction of
  spaces with given homotopy groups, via the Dold-Kan correspondence
  (compare~\cite{ELST16:trivial}).
\end{remark}

\section{Applications}\label{sec:apps}
In addition to the detection results described in Section~\ref{sec:ssec}, several of the other most celebrated
applications of Khovanov homology also come from the spectral
sequences discussed in Section~\ref{sec:ssec}, though there are other
important applications not directly tied to these spectral
sequences. Like the rest of the paper, our intention in this section
is to give a sense of the breadth of applications of these techniques,
and some of the ideas behind them, not a comprehensive list.

Given a knot $K$, Rasmussen observed that the two copies of $\Q$ in the $E^\infty$-page of the Lee spectral sequence lie in adjacent quantum gradings $s(K)\pm 1$. (Recall that the quantum gradings for a knot are always odd.) He further showed that:
\begin{theorem}(Rasmussen \cite{Rasmussen10:s})
  The integer $s(K)$ is a homomorphism from the smooth concordance
  group onto $2\Z$. Further, if there is a genus $g$ knot cobordism
  from $K$ to $K'$ then $|s(K)-s(K')|\leq 2g$.
\end{theorem}
The proof is combinatorial and relatively simple: the K\"unneth theorem for Lee homology implies that the Rasmussen invariant $s(K)$ is additive for connected sums, and the fact that the two copies of $\Q$ are in adjacent gradings quickly gives that $s(m(K))=-s(K)$ (where $m$ denotes the mirror). Rasmussen then shows that the maps on Khovanov homology associated to elementary cobordisms $\Sigma$ lift to maps of the Lee complex changing the filtration by $-\chi(\Sigma)$, and that the map associated to a connected cobordism is an isomorphism on Lee homology. The result follows.

Rasmussen's construction was inspired by Ozsv\'ath-Szab\'o's $\tau$ invariant~\cite{OS03:4BallGenus}. In fact, Rasmussen initially conjectured that $s(K)$ was equal to $2\tau(K)$, but this conjecture was quickly disproved, showing that, in fact, $s(K)$ and $\tau(K)$ together gave the first surjection from the smooth concordance group of topologically slice knots to $\Z^2$~\cite{HO08:tau-not-s,Livingston08:s-tau}. Similar constructions have been given using other spectral sequences, including several infinite families of concordance invariants~\cite{SSS17:perturb,LL19:upsilonish}, though these are not known to be independent. (Again, these families were inspired by constructions in Heegaard Floer theory which, in that case, were shown to give a surjection from the smooth concordance group of topologically slice knots onto $\mathbb{Z}^\infty$~\cite{OSS17:upsilon,DHST:concordance}.)

For certain classes of knots, the $s$-invariant is easy to compute. In particular, this holds for positive knots, i.e., knots where all the crossings are positive. As a corollary, Rasmussen obtains the following remarkable generalization of the Milnor conjecture on the slice genus of torus knots:
\begin{corollary}(Rasmussen \cite{Rasmussen10:s})
  Let $K$ be a positive knot with $n$ crossings and where the oriented resolution of $K$ has $k$ circles. Let $g_4(K)$ be the slice genus of $K$ and $g(K)$ the ordinary knot genus of $K$. Then
  \[
    s(K)=2g_4(K)=2g(K)=n-k+1.
  \]
\end{corollary}

For the case of torus knots, the equality $g_4(T_{p,q})=g(T_{p,q})=\frac{(p-1)(q-1)}{2}$ was conjectured by Milnor~\cite{Milnor68:singular} and first proved by Kronheimer-Mrowka, by applying gauge theory to bound the genera of embedded surfaces in the K3-surface~\cite{KM93:Milnor1,KM95:Milnor2}. It also follows from Thom's conjecture about the genera of embedded surfaces in $\C P^2$, first proved using Seiberg-Witten gauge theory~\cite{KM94:Thom,MST96:Thom}. Rasmussen's argument is the first combinatorial proof of the Milnor conjecture. As he observes, the $s$ invariant is not a lower bound on the topological slice genus, and in fact can be used to show that some topologically slice knots are not smoothly slice. (See~\cite{Plam06:transverse,Shumakovitch07:s} for conceptual proofs of this fact.) When combined with work for M.~Freedman~\cite{Freedman82,FQ90}, this also implies the existence of exotic smooth structures on $\R^4$.

Remarkably, Lambert-Cole recently showed that one can deduce the Thom conjecture from Milnor's conjecture~\cite{LC20:Thom}, so Rasmussen's result underlies a combinatorial proof of the Thom conjecture as well. (The proof uses and extends Gay-Kirby's notion of \emph{trisections}~\cite{GK16:trisections}.) Pushing these ideas further, he was even able to recover the generalized Thom conjecture, that pseudo-holomorphic curves in symplectic 4-manifolds---and, in particular, in K\"ahler manifolds like $\C P^2$---are genus-minimizing~\cite{LC:symp-Thom} (again originally proved using gauge theory~\cite{OS00:symp-Thom}).

Another striking application of the $s$-invariant was given recently by L.~Piccirillo, who used it to show that the Conway mutant of the Kinoshita-Terasaka knot is not smoothly slice, resolving a longstanding question~\cite{Piccirillo20:Conway}. (By work of Freedman, any knot with Alexander polynomial $1$ is topologically slice~\cite{FQ90}, and the Kinoshita-Terasaka knot is smoothly slice. Indeed, the Conway knot was the only knot with 13 or fewer crossings whose slice status was not known.) Piccirillo recalls that a knot $K$ is smoothly slice if and only if the $0$-trace of $K$, the result of attaching a $0$-framed $2$-handle to the 4-ball along $K$, embeds smoothly in $S^4$.  The $s$-invariant of the Conway knot vanishes, but Piccirillo produces another knot $K'$ whose zero trace is diffeomorphic to the $0$-trace of the Conway knot, as can be shown by explicit handle calculus such that $s(K')\neq0$. As she notes, the $s$-invariant plays a special role here: other known smooth concordance invariants, like the Heegaard Floer analogue $\tau$, would not work for this strategy. This proof, and the $s$-invariant, gives a possible attack on the smooth 4-dimensional Poincar\'e conjecture~\cite{MP:exotic-s} (see also~\cite{FGMW10:man-and-machine,MMSW:S1S2}).

Functoriality of Khovanov homology means that it also gives an invariant of surfaces in $\mathbb{R}^4$. For closed surfaces, this invariant turns out not to be interesting: it vanishes if some component of the surface is not a torus, and otherwise is $2^n$ if the surface consists of $n$ tori~\cite{Rasmussen05:kh-closed,Tanaka06:closed,LG:split}. On the other hand, for surfaces with boundary a nontrivial link in $S^3$, Khovanov homology does give an interesting invariant~\cite{SS:relative}, even distinguishing some surfaces that are topologically isotopic~\cite{HS:surfaces,LS:mixed}.

In a different direction, Khovanov homology and its cousins have had interesting applications to Legendrian and transverse knot theory. Recall that a knot $K(t)=(x(t),y(t),z(t))$ in $\R^3$ is \emph{Legendrian} if $y(t)=z'(t)/x'(t)=dz/dx$ for all $t$, and \emph{transverse} if this condition holds for no $t$. A Legendrian knot $K$ has three \emph{classical invariants}: its underlying smooth knot type; its Thurston-Bennequin number $\mathrm{TB}(K)$, which is the difference between the Seifert framing and the framing induced by the 2-plane field $\ker(dz-ydx)$; and the \emph{rotation number} $\mathrm{rot}(K)$, which is the relative Euler class of the 2-plane field over a Seifert surface.
(See~\cite{Etnyre05:survey} for a nice survey on Legendrian and transverse knots.) Given a Legendrian knot $K$, there are stabilization operations that do not change the underlying smooth knot but change the pair $(\mathrm{TB}(K),\mathrm{rot}(K))$ by $(-1,\pm 1)$.
The celebrated \emph{slice-Bennequin} inequality states that for a given smooth knot type, $\mathrm{TB}(K)+|\mathrm{rot}(K)|\leq 2g_4(K)-1$, where $g_4(K)$ is the slice genus~\cite{Bennequin83,Rudolph93:slice-Ben}. So, the pairs $(\mathrm{TB}(K),\mathrm{rot}(K))$ realized by Legendrian representatives of a smooth knot type form a mountain range.

L.~Ng improved the slice-Bennequin inequality to show that
\begin{equation}\label{eq:Ng-bound}
  \min\Bigl\{k\mid \bigoplus_{i-j=k}\Kh^{i,j}(K)\neq 0\Bigr\}
\end{equation}
is an upper bound on the Thurston-Bennequin number of any Legendrian representative of $K$~\cite{Ng05:TB1} (see also\cite{Shumakovitch07:s}). In particular, this gives the bound $\mathrm{TB}(K)\leq s(K)-1$, a refinement of the slice-Bennequin inequality. The bound~\eqref{eq:Ng-bound} is sharp for alternating knots, and by combining it with tools from Heegaard Floer homology Ng computed the maximal Thurston-Bennequin number for all knots up to 11 crossings~\cite{Ng12:TB2}.
In fact, many different bounds on the Thurston-Bennequin number, including this one and another coming from the Kauffman polynomial, have a common skein-theoretic proof~\cite{Ng08:TB3}. (Both the Kauffman polynomial and Khovanov homology bounds are often sharp, at least for small knots.)
As another potential application to contact topology, O.~Plamenevskaya defined a natural invariant of transverse knots (and, consequently, Legendrian knots) lying in Khovanov homology~\cite{Plam06:transverse}. Several variants of her construction have been given~\cite{LNS15:transverse,Wu08:transverse,HS16:transverse,Montes}, but it remains open whether any of these invariants is effective, i.e., distinguishes some pair of transverse knots with the same classical invariants.

A third class of application has been to ribbon cobordisms. Generalizing the notion of a ribbon knot, C.~Gordon introduced the notion of a ribbon concordance from a knot $K_1$ to a knot $K_2$: a concordance is \emph{ribbon} if it is built entirely from births and saddles~\cite{Gordon81:ribbon}. So, a knot $K$ is ribbon if there is a ribbon concordance from the unknot to $K$. Ribbon concordance is not symmetric, and a conjecture of Gordon's, recently proved by Agol~\cite{Agol:Gordon-true}, is that ribbon concordance forms a partial order: if $K_1$ is ribbon concordant to $K_2$ and $K_2$ is ribbon concordant to $K_1$ then $K_1=K_2$. Inspired by an analogous result for Heegaard Floer homology~\cite{Zemke19:ribbon}, A.~Levine and I.~Zemke show that a ribbon concordance induces a split injection on Khovanov homology~\cite{LZ19:ribbon}. In particular, this implies that if there is a ribbon concordance from an alternating knot $K_1$ to $K_2$ then the crossing number of $K_2$ is at least as large as the crossing number of $K_1$---an elementary obstruction of which no other proof is currently known.

For many of these applications to be effective, one needs an efficient way to compute Khovanov homology and, ideally, the Lee spectral sequence. There are several programs that compute versions of Khovanov homology directly~\cite{KhoHo,KnotKit,KnotJob}. Since the Khovanov cube itself grows exponentially, direct computations become impossible around 17 crossings. Fortunately, the tangle invariants provide more efficient algorithms, through an approach that Bar-Natan calls \emph{scanning}~\cite{BarNatan07:fast,FastKh}, an idea that, on the decategorified level, goes back to  Jones and his work on the Temperley-Lieb algebra~\cite{Jones91:subfactors-book}. First, you factor a knot as $T_1T_2\cdots T_k$, and compute the invariant of each $T_i$. You tensor the invariants for $T_1$ and $T_2$, simplify the result, then tensor on the invariant for $T_3$, simplify the result, and so on. This allows one to compute the invariant for much larger knots, including the $s$-invariant for a 78-ish crossing knot of interest~\cite{FGMW10:man-and-machine} and the 49-crossing knot needed for Piccirillo's proof described above.

\vspace{0.1in}

\subsection*{Acknowledgments} We thank Mohammed Abouzaid, John Baldwin, Joshua Greene, Ciprian Manolescu, and Sucharit Sarkar for helpful comments. MK was partially supported by NSF grant DMS-1807425 while working on this paper; RL was partially supported by NSF grant DMS-1810893.


\bibliographystyle{hamsplain}\bibliography{jones}

\end{document}